\documentclass{article}
\usepackage[utf8]{inputenc}
\usepackage{graphicx}
\usepackage{amsmath}
\usepackage{hyperref}
\usepackage{tikz}
\usepackage{geometry}
\usepackage{multicol}
\usepackage[numbers]{natbib}
\usepackage{amssymb,amsfonts,fullpage,amsthm,color,mathrsfs, authblk}

\theoremstyle{definition}
\newtheorem{Theorem}{Theorem}

\newtheorem{Lemma}{Lemma}
\numberwithin{Lemma}{section}
\numberwithin{Theorem}{section}
\numberwithin{Proposition}{section}
\numberwithin{Definition}{section}
\numberwithin{Remark}{section}

\newcommand{\R}{\mathbb{R}}

\newcommand{\N}{\mathbb{N}}
\newcommand{\C}{\mathbb{C}}

\newcommand{\D}{\mathcal{D}}
\newcommand{\so}{\mathfrak{so}}
\newcommand{\Uq}{\mathcal{U}_q}

\makeatletter
\def\blfootnote{\xdef\@thefnmark{}\@footnotetext}
\makeatother

\title{A Type D Asymmetric Simple Exclusion Process Generated by an Explicit Central Element of $\Uq(\so_{10})$}
\date{}

\author[1]{Eddie Rohr}
\author[2]{Karthik Sellakumaran Latha}
\author[3 *]{Amanda Yin}
\affil[1]{Undergraduate, Department of Mathematics, The College of New Jersey}
\affil[2]{Undergraduate, Department of Mathematics, University of Maryland}
\affil[3]{Undergraduate, Department of Mathematics, University of Texas at Austin}

\begin{document}
\maketitle
\blfootnote{*Author emails: rohre1@tcnj.edu, karthiks@umd.edu, amandayin@utexas.edu}

\vspace{-1cm}
\begin{abstract}
    The Type D asymmetric simple exclusion process (Type D ASEP) is a two-species interacting particle system exhibiting a drift, where two particles may occupy the same site only if they belong to different species. In previous research \cite{kuan2020}, the Type D ASEP was generated using the quantum Hamiltonian corresponding to central elements from the quantum groups $\Uq(\so_6)$ and $\Uq(\so_8)$. We extend this construction to the case of $\Uq(\so_{10})$.  Additionally, we generalize a previously known duality function from \cite{Blyschak_2023} for the Type D ASEP for all $n$.
\end{abstract}

\section{Introduction}

The \textbf{Type D asymmetric simple exclusion process (Type D ASEP)} is an interacting particle system introduced in \cite{kuan2020} which generalizes the asymmetric simple exclusion process originally introduced by Spitzer \cite{Spit70}. The state space of the Type D ASEP consists of two species of particles interacting on a one-dimensional lattice, where two particles may occupy the same site as long as they are of different species. The Type D ASEP has three parameters $(q, n, \delta)$, where $q$ is the asymmetry parameter, the positive integer $n$ gives the speed of the drift, and $\delta$ specifies the interaction between the two particle species. The generator matrix is given in Section \ref{subsec:asep}.

The procedure connecting quantum groups to a particle system and duality is described in \cite{CGRS}; see also \cite{Sch97, IS11, CGRS16a}. Previous research \cite{Blyschak_2023, kuan2020} uses this procedure to construct Type D ASEPs with parameters $n = 3, 4$ from the representation theory of the Type D quantum groups $\Uq(\so_6)$ and $ \Uq(\so_8)$, respectively. The authors of \cite{kuan2020} conjecture that this process can be generalized to $\Uq(\so_{2n})$. 

This paper verifies this claim in the case of $\Uq(\so_{10})$. The first step toward this is to compute a central element of this quantum group using methods from \cite{kuan2020}. Then, a generator of the Type D ASEP with parameter $n = 5$ is constructed from this central element. 

In addition, this paper expands on previous results concerning Markov self-duality of the Type D ASEP. In \cite{kuan2020}, the authors produce a triangular duality function. This paper instead focuses on orthogonal polynomial duality; see \cite{Gro19, CFG19, Franceschini_2018, franceschini2022orthogonal, carinci2021qorthogonal}. In particular, \cite{Blyschak_2023} presents a duality function for Type D ASEP with parameters $n = 3, 4$ and $\delta = 0$. This paper generalizes this result for all $n$. 

The paper is outlined as follows: Section \ref{sec:background} describes the relevant background and notation used. Section \ref{sec:results} states the results: a central element of $\Uq(\so_{10})$, the construction of the Type D ASEP generator for $n = 5$, and a self-duality function of the Type D ASEP. Section \ref{sec:proofs} outlines the proofs of these results.

\textbf{Acknowledgments.}
This research was conducted as a part of the Texas A\&M University math REU, funded in part by the National Science Foundation (DMS--2150094). The authors would like to thank their mentor Dr. Jeffrey Kuan, their teaching assistant Zhengye Zhou, and Andrew Lin for their support and involvement.

\section{Background}\label{sec:background}

\subsection{Algebraic Definitions}
The \textbf{special orthogonal Lie algebra} $\so_{2n}$ is defined as
$$\so_{2n} = \Biggl\{\begin{pmatrix}
A & B \\
C & -A^T 
\end{pmatrix}
\Bigg|A,B,C \in M_{n\times n}(\C), B = -B^T, C = -C^T\Biggl\}.$$

Let $E_{i,j}$ denote the matrix with $1$ in entry $(i, j)$ and $0$ elsewhere. Then define matrices $E_i, F_i, H_i, 1 \leq i \leq n$ as follows:
\begin{align*}
    E_n &= E_{n-1, 2n} - E_{n, 2n-1} \\
    F_n &= E_{2n-1, n} - E_{2n, n-1} \\
    H_n &= E_{n-1, n-1} + E_{n, n} - E_{2n - 1, 2n-1} - E_{2n,2n} 
\end{align*}
and for all $1 \leq i \leq n - 1$,
\begin{align*}
    E_i &= E_{i, i+1} - E_{n+i+1, n+i} \\
    F_i &= E_{i+1, i} - E_{n+i, n+i+1} \\
    H_i &= E_{i,i} - E_{i+1, i+1} - E_{n+i,n+i} + E_{n+i+1,n+i+1}.
\end{align*}
These matrices generate the fundamental representation of $\so_{2n}$.

This Lie algebra has rank $n$ and Dynkin diagram $D_n$, with Cartan matrix $(a_{ij})$ given by
$$a_{ij} = \begin{cases}
    2 & i = j \\
    -1 & \{i, j\} = \{n - 2, n\} \text{ or } \{k, k+1\}, 1 \leq k \leq n - 2 \\
    0 & \text{otherwise}
\end{cases}.$$

Finally, let $L_i$ denote the linear operator taking a matrix to its $i$-th diagonal entry. Then the positive roots of $\so_{2n}$ are $\{L_i+L_j\}_{i<j}\cup \{L_i-L_j\}_{i<j}$. The simple roots are $\alpha_i = L_i - L_{i+1}$ for $1 \leq i \leq n - 1$ and $\alpha_n = L_{n-1} + L_n$. In the fundamental representation of $\so_{2n}$, elements of the Lie algebra act as the underlying matrix on $\C^{2n}$. The fundamental weights are then $\pm L_i$, $1 \leq i \leq n$.

An ordering for the weights 
\[
    L_1 > L_2 >\ldots > L_n = -L_n > -L_{n-1} >\ldots > -L_1
\] is given in \cite{kuan2020}.
\\~\

Define vectors $v_1, \ldots, v_{2n}$ with $v_i$ having a $1$ in the $i$-th coordinate if $1 \le i \le n$, a $1$ in the $(3n + 1 - i)$-th coordinate if $n+1 \le i \le 2n$, and zeroes elsewhere. Observe that $v_1, \ldots, v_{2n}$ are in the weight spaces of $L_1, \ldots, L_n, -L_n, \ldots, -L_1$, respectively.

The \textbf{quantum group} $\Uq(\so_{2n})$ is the algebra generated by $\{E_i, F_i, q^{H_i} : 1\leq i \leq n\}$ with $q$-deformed relations
$$
[E_i, F_i] = \tfrac{q^{H_i} - q^{-H_i}}{q-q^{-1}},\quad q^{H_i}E_j = q^{(\alpha_i,\alpha_j)}E_jq^{H_i},\quad q^{H_i}F_j = q^{-(\alpha_i,\alpha_j)}F_jq^{H_i},
$$
along with the Serre relation for every $(i,j)$ with $a_{ij} = -1$ in the Cartan matrix
$$
E_i^2E_j + E_jE_i^2 = (1+q)E_iE_jE_i,\quad F_i^2F_j + F_jF_i^2 = (1+q)F_iF_jF_i.
$$
Here, $(\alpha_i,\alpha_j)$ is the usual dot product.
All other pairs of elements commute.
Lastly, the coproducts of the generators are 
\[
    \Delta (E_i) = E_i \otimes 1 + q^{H_i} \otimes E_i, \quad \Delta(F_i) = 1 \otimes F_i + F_i \otimes q^{-H_i}, \quad \Delta(q^{H_i}) = q^{H_i} \otimes q^{H_i}.
\]
For convenience, in the remainder of the paper, let $K_i$ denote $q^{H_i}$.

A general definition for the quantum group of a Lie algebra $\mathfrak{g}$ is given in \cite{Jan95}, and in Chapter 6 the author introduces a bilinear pairing that will be used in the later construction. Recall that the Borel subalgebras $\mathfrak{b}\pm$ are the Lie subalgebras generated by $\{E_i, H_i\}$ and $\{F_i, H_i\}$, respectively. Similarly, let $\Uq(\mathfrak{b}\pm)$ denote the corresponding subalgebras of the quantum group generated by the Borel subalgebras (replacing $H_i$ with $q^{H_i}$). There is a bilinear pairing $\langle\cdot \, , \,\cdot \rangle$: $\Uq(\mathfrak{b}-) \times \Uq(\mathfrak{b}+) \to \mathbb{Q}(q)$ such that for any linear combinations $\alpha, \beta$ of the positive simple roots we have
\[
    \langle q^{H_\alpha}, q^{H_\beta} \rangle = q^{-(\alpha \cdot \beta)} \text{ and } \langle F_i, E_j \rangle = \frac{-\delta_{ij}}{(q-\frac{1}{q})},
\]
and all other pairings between generators are zero. The pairing can be computed for products via
\[
    \langle y, xx' \rangle = \langle \Delta(y), x' \otimes x \rangle, \quad \langle yy', x \rangle = \langle y \otimes y', \Delta(x) \rangle,
\]
where $\Delta(xy) = \Delta(x)\Delta(y)$, and $\langle x_1 \otimes x_2, y_1 \otimes y_2 \rangle$ is defined to be $\langle x_1, y_1 \rangle \langle x_2, y_2 \rangle$.

\subsection{Probabilistic definitions}\label{subsec:asep}

Recall that the Type D ASEP consists of two species of particles interacting on a one-dimensional lattice, where at most two particles can occupy a single site, and if so must be of different species. The Type D ASEP has three parameters $(q, n, \delta)$, where $q \neq 1$ is a positive real number, and $n\geq 2$ and $0\leq \delta \leq n-2$ are integers. Intuitively, $q$ is the asymmetry parameter, $n$ affects the speed of the drift, and $\delta$ affects the interaction between the two species. If the lattice has $L$ sites, then there are $4^L$ states. A more complete definition of the Type D ASEP is given in \cite{kuan2020}.\footnote{Dr. Kuan and Zhengye, two authors of \cite{Blyschak_2023, kuan2020}, acknowledged a typo in those papers where $\so_{2n}$ was associated with the parameters $(q,n-1,\delta)$ rather than $(q,n,\delta)$.}

First, consider the two-site model. The $16$ states are $(a, b)$, $0 \leq a, b \leq 3$, where each coordinate corresponds to a site, $0$ denotes an empty site, $1$ denotes a species $1$ particle, $2$ denotes a species $2$ particle, and $3$ denotes both a species $1$ and species $2$ particle. The generator matrix is given by the direct sum of a $4 \times 4$ block   
\renewcommand*{\arraystretch}{2.5}
\begin{equation*}
\mathcal{L}_1 = \begin{small}\begin{bmatrix}* & \frac{\left(1 + q^{- 2 \delta}\right) \left(q^{2 n - 2} - q^{2 n -4} + \frac{2}{q^{2}}\right)}{q^{2 \delta} + 1} & \frac{\left(- q^{1 - n} + q^{n - 1}\right)^{2}}{q^{2}} & \frac{q^{2 \delta} \left(1 + q^{- 2 \delta}\right) \left(q^{2n - 2} - q^{2 n - 4} + \frac{2}{q^{2}}\right)}{q^{2 \delta} + 1}\\
q^{- 2 \delta} \left(q^{2 n} - q^{2 n - 2} + 2\right) & * & q^{-2 n} - q^{2-2 n} + 2 & \left(- q^{1 - n} + q^{n - 1}\right)^{2}\\
q^{2} \left(- q^{1 - n} + q^{n - 1}\right)^{2} & 2 q^{2} + q^{2 - 2 n} - q^{4 - 2 n} & * & q^{2 \delta} \left(2 q^{2} + q^{2 - 2 n} - q^{4 - 2 n}\right)\\
q^{2 n} - q^{2 n - 2} + 2 & \left(- q^{1 - n} + q^{n - 1}\right)^{2} & q^{2 \delta} \left(q^{-2 n} - q^{2 - 2 n} + 2\right) & *\end{bmatrix}\end{small}
\end{equation*}
\renewcommand*{\arraystretch}{1.5}
    corresponding to the communicating class $\{(3, 0), (2, 1), (0, 3), (1, 2)\}$, four $2 \times 2$ blocks
    \begin{equation*}
        \mathcal{L}_2 = \begin{bmatrix}
        * & \tfrac{q^{1 - 2n} + q^{2n - 1}}{q} \\
        q(q^{1 - 2n} + q^{2n - 1}) & *
    \end{bmatrix}
    \end{equation*}
    corresponding to the communicating classes $\{(1, 0), (0,1)\}, \{(2, 0), (0, 2)\}, \{(3, 1), (1, 3)\}, \{(3, 2), (2, 3)\}$, and four $1 \times 1$ blocks with entry $0$ corresponding to the communicating classes $\{(0, 0)\}, \{(1, 1)\}, \{(2, 2)\}, \{(3, 3)\}$. Here, the diagonal entries are chosen so that the rows sum to $0$. To summarize, the generator matrix is
    \begin{equation}\label{eq:generator}
        \mathcal{L} = \mathcal{L}_1 \oplus \bigoplus_{i=1}^4 \mathcal{L}_2 \oplus \bigoplus_{i=1}^4 [0],
    \end{equation}
    with respect to the ordered basis $\{(3, 0), (2, 1), (0, 3), (1, 2), (1, 0), (0,1), (2, 0), (0, 2), (3, 1), (1, 3), (3, 2),\\ (2, 3), (0, 0), (1, 1), (2, 2), (3, 3)\}$.
    
    If there are $L$ sites, then the generator matrix is given by
    $$\mathcal{L}^{1,2} + \mathcal{L}^{2,3} + \dots + \mathcal{L}^{L - 1, L},$$
    where $\mathcal{L}^{x,x+1}$ denotes the matrix acting on lattice sites $x$ and $x + 1$. 

    We now introduce several notations and definitions that will appear in Theorem \ref{thm:duality}. If $x \in \{1, \dots, L\}$ is a lattice site, $i \in 1, 2$, and $\eta$ is a state, then let $\eta_i^x$ denote the number of particles of class $i$ at site $x$ in state $\eta$ (so $\eta_i^x$ is $0$ or $1$). Define height functions
    $$N^-_{x-1} (\eta_i) = \sum_{i=1}^{x-1} \eta_i^x, \quad N^+_{x+1} (\eta_i) = \sum_{i=x+1}^{L} \eta_i^x$$
    that count the number of particles to the left and right of a particular site.

    Define the $q$-Pochhammer symbol for $a \in \R, k \in \N$ as
    $$(a; q)_k := \prod_{i=0}^{k-1} (1 - aq^i),$$
    define the $q$-hypergeometric function $_2\varphi_1$ as
    $$_2\varphi_1\left(\begin{matrix}a,b\\c\end{matrix}; q, z\right) := \sum_{k=0}^\infty \frac{(a;q)_k (b;q)_k}{(c;q)_k} \frac{z^k}{(q;q)_k},$$
    and define the $q$-Krawtchouk polynomials as
    $$K_n(q^{-x}; p, c; q) = {}_2\varphi_1 \left(\begin{matrix}q^{-x},q^{-n}\\q^{-c}\end{matrix}; q, pq^{n+1}\right).$$

\section{Results}\label{sec:results}
First, for ease of reading set some notational shortcuts following the convention in \cite{kuan2020}. Set $r = q - \frac{1}{q}$, and let $E_{x_1x_2\ldots x_k} := E_{x_1}E_{x_2}\cdots E_{x_k}$, and similarly for $F$ and $K$. For example, $r^2F_{23}E_{23} = (q-\frac{1}{q})^2F_2F_3E_2E_3.$
\begin{Theorem}\label{thm:central-element}
    The following element is in the center of $\Uq(\so_{10})$:
    \begin{align*}
   &q^8  K_{11223345} + q^6 K_{223345} + q^4 K_{3345} + q^2 K_{45} + K_{4} K_5^{-1} + K_4^{-1} K_{5} + \tfrac{1}{q^2} K_{45}^{-1} + \tfrac{1}{q^4} K_{3345}^{-1} + \tfrac{1}{q^6} K_{223345}^{-1} \\
   &+ \tfrac{1}{q^8} K_{11223345}^{-1} +\tfrac{r^2}{q} F_{4} K_5^{-1} E_{4} + \tfrac{r^2}{q^3} (q F_{34} - F_{43}) K_{35}^{-1} (q E_{43} - E_{34}) + \tfrac{r^2}{q^3} F_{3} K_{345}^{-1} E_{3} \\
   &+ \tfrac{r^2}{q^5} (q^2 F_{234} - q F_{243} - q F_{324} + F_{432}) K_{235}^{-1} (q^2 E_{432} - q E_{342} - q E_{423} + E_{234})\\
    &+ \tfrac{r^2}{q^5} ({q} F_{23} - F_{32}) K_{2345}^{-1} ({q} E_{32} - E_{23}) + \tfrac{r^2}{q^5} F_{2} K_{23345}^{-1} E_{2} + \tfrac{r^2}{q^7} (\boxed{A_1}) K_{1235}^{-1} (\boxed{A_2}) \\
    &+ \tfrac{r^2}{q^7} ({q^2} F_{123} - {q} F_{132} - {q} F_{213} + F_{321})  K_{12345}^{-1} ({q^2} E_{321} - {q} E_{231} - {q} E_{312} + E_{123})\\
    &+ \tfrac{r^2}{q^7} ({q} F_{12} - F_{21})  K_{123345}^{-1}  ({q} E_{21} - E_{12}) + \tfrac{r^2}{q^7} F_{1} K_{1223345}^{-1} E_{1} - q r^{2} (\boxed{A_3}) K_{1234} (\boxed{A_4}) \\
        &- q r^{2} (q^2 F_{532} - q F_{352} - q F_{523} + F_{235}) K_{234} (q^2 E_{235} - q E_{253} - q E_{325} + E_{532})\\ 
        &- q r^{2} (q F_{53} - F_{35}) K_{34} (q E_{35} - E_{53}) - q r^{2} F_{5}    K_{4} E_{5} - r^{4} F_{54} E_{54} - \tfrac{r^2}{q} (\boxed{A_5}) K_{123} (\boxed{A_6}) - \tfrac{r^2}{q} (\boxed{A_7}) K_{23} (\boxed{A_8})\\
        &- \tfrac{r^{2}}{q} (q^2 F_{453} - q F_{435} - q F_{534} + F_{345}) K_{3} (q^2 E_{354} - q E_{435} - q E_{534} + E_{543}) - \tfrac{r^{2}}{q} F_{5} K_{4}^{-1} E_{5} \\
        &- \tfrac{r^4}{q^2} (- q F_{3435} - q F_{5343} + (q^2 + 1) F_{3543}) (- q E_{3435} - q E_{5343} + (q^2 + 1) E_{3543})- \tfrac{r^{2}}{q^3} (\boxed{A_9}) K_{2} (\boxed{A_{10}})\\ 
        &- \tfrac{r^2}{q^3} (q^2 F_{354} - q F_{435} - q F_{534} + F_{543}) K_{3}^{-1} (q^2 E_{453} - q E_{435} - q E_{534} + E_{345}) - \tfrac{r^2}{q^3} (q F_{35} - F_{53}) K_{34}^{-1} (q E_{53} - E_{35})\\ 
        &- \tfrac{r^2}{q^3} (\boxed{A_{11}}) K_{12} (\boxed{A_{12}}) - \tfrac{r^{4}}{q^4} (\boxed{A_{13}}) (\boxed{A_{14}}) - \tfrac{r^{2}}{q^5} (\boxed{A_{15}}) K_{2}^{-1} (\boxed{A_{16}}) - \tfrac{r^2}{q^5} (\boxed{A_{17}}) K_{23}^{-1} (\boxed{A_{18}})\\
        &- \tfrac{r^{2}}{q^5} (q^2 F_{235} - q F_{253} - q F_{325} + F_{532}) K_{234}^{-1} (q^2 E_{532} - q E_{352} - q E_{523} + E_{235}) - \tfrac{r^{2}}{q^5} (\boxed{A_{19}}) K_{1} (\boxed{A_{20}}) \\
        & - \tfrac{r^{4}}{q^6} (\boxed{A_{21}}) (\boxed{A_{22}}) - \tfrac{r^{2}}{q^7} (\boxed{A_{23}}) K_{123}^{-1} (\boxed{A_{24}}) - \tfrac{r^2}{q^7} (\boxed{A_{25}}) K_{1234}^{-1} (\boxed{A_{26}}) - \tfrac{r^{2}}{q^7} (\boxed{A_{27}}) K_{12}^{-1} (\boxed{A_{28}}) \\ 
        &- \tfrac{r^{2}}{q^7} (\boxed{A_{29}}) K_{1}^{-1} (\boxed{A_{30}}) q^7 r^2 F_{1} K_{1223345} E_{1} + q^5 r^2 ({q} F_{21} - F_{12}) K_{123345} ({q} E_{12} - E_{21}) + q^5 r^2 F_{2}  K_{23345} E_{2}\\
        &+ q^3 r^2 ({q^2} F_{321} - {q} F_{231} - {q} F_{312} + F_{123})  K_{12345} ({q^2} E_{123} - {q} E_{132} - {q} E_{213} + E_{321})\\
        &+ q^3 r^{2} ({q} F_{32} - F_{23}) K_{2345} ({q} E_{23} - E_{32}) + q^3 r^2 F_{3}  K_{345} E_{3} + q r^2 (\boxed{A_{31}}) K_{1235} (\boxed{A_{32}})\\
        &+ q r^2 ({q^2} F_{432} - {q} F_{342} - {q} F_{423} + F_{234}) K_{235} ({q^2} E_{234} - {q} E_{243} - {q} E_{324} + E_{432}) \\
        &+ q r^2 ({q} F_{43} - F_{34}) K_{35} ({q} E_{34} - E_{43}) + q r^2 F_{4} K_{5} E_{4},
\end{align*}
where
\begin{align*}
    A_1 &= {q^3} F_{1234} - {q^2} F_{1243} - {q^2} F_{1324} - {q^2} F_{2134} + {q} F_{1432} + {q} F_{2143} + {q} F_{3214} - F_{4321} \\
    A_2 &= {q^3} E_{4321} - {q^2} E_{3421} - {q^2} E_{4231} - {q^2} E_{4312} + {q} E_{2341} + {q} E_{3412} + {q} E_{4123} - E_{1234} \\
    A_3 &= q^3 F_{5321} - q^2 F_{3521} - q^2 F_{5231} - q^2 F_{5312} + q F_{2351} + q F_{3512} + q F_{5123} - F_{1235} \\
    A_4 &= q^3 E_{1235} - q^2 E_{1253} - q^2 E_{1325} - q^2 E_{2135} + q E_{1532} + q E_{2153} + q E_{3215} - E_{5321} \\
    A_5 &= q^4 F_{45321} - q^3 F_{43521} - q^3 F_{45231} - q^3 F_{45312} - q^3 F_{53421} + q^2 F_{34521} + q^2 F_{42351} + q^2 F_{43512}\\
    &    + q^2 F_{45123} + q^2 F_{52341} + q^2 F_{53412} - q F_{23451} - q F_{34512} - q F_{41235} - q F_{51234} + F_{12345} \\
    A_6 &= q^4 E_{12354} - q^3 E_{12435} - q^3 E_{12534} -q^3 E_{13254} - q^3 E_{21354} + q^2 E_{12543} + q^2 E_{14325} + q^2 E_{15324} \\
    &    + q^2 E_{21435} + {q^2} E_{21534} + q^2 E_{32154} - q E_{15432} - q E_{21543} - q E_{43215} - q E_{53214} + E_{54321} \\
    A_7 &= q^3 F_{4532} - q^2 F_{4352} - q^2 F_{4523} - q^2 F_{5342} + q F_{3452} + q F_{4235} + q F_{5234} - F_{2345} \\
    A_8 &= q^3 E_{2354} - q^2 E_{2435} - q^2 E_{2534} - q^2 E_{3254} + q E_{2543} + q E_{4325} + q E_{5324} - E_{5432} \\
    A_9 &= q^4 F_{34532} - (q^3 - q) F_{35342} - (q^3 - q) F_{43532} + q^2 F_{34235} + q^2 F_{35234} + q^2 F_{43523}\\
    &    + q^2 F_{53423} - q F_{32345} - q F_{43235} - q F_{45323} - q F_{53234} + F_{23453} + (- q^3 - q) F_{34523} \\
    A_{10} &= q^4 E_{23543} - q^3 E_{23435} - q^3 E_{25343} + q^2 E_{32435} + q^2 E_{32534} + q^2 E_{43253}+ q^2 E_{53243} \\
    &     - q E_{34325} - q E_{53432} + E_{35432} + (- q^3 - q) E_{32543} 
\end{align*}
and
\begin{align*}
    A_{11} &= q^5 F_{354321} - q^4 F_{354312} - (q^4 - q^2 ) F_{343521} - (q^4 - q^2 ) F_{534321} + q^3 F_{342351} + q^3 F_{352341}  \\
    &    + q^3 F_{435231} + q^3 F_{534231} - q^2 F_{323541} - q^2 F_{341235} - q^2 F_{351234} - q^2 F_{432351} - q^2 F_{435123} \\
    &    - q^2 F_{532341} - q^2 F_{534123} - q^2 F_{543231} + q F_{235431} + q F_{312354} + q F_{431235} + q F_{531234} \\
    &    + q F_{543123} + (q^3 - q) F_{343512} + (q^3 - q) F_{534312} + (q^3 + q) F_{354123} + (- q^4 - q^2) F_{354231} - F_{123543} \\
    A_{12} &= q^5 E_{123543} - q^4 E_{123435} - q^4 E_{125343} - q^4 E_{213543} + q^3 E_{132435} + q^3 E_{132534} + q^3 E_{143253}\\
    &     + q^3 E_{153243} + q^3 E_{213435} + q^3 E_{215343} - q^2 E_{134325} - q^2 E_{153432} - q^2 E_{321435} - q^2 E_{321534} - q^2 E_{432153}\\
    &      - q^2 E_{532143} + q E_{135432} + q E_{343215}   + q E_{534321} + (q^3 + q) E_{321543} + (- q^4 - q^2) E_{132543} - E_{354321} \\
    A_{13} &= q^2 F_{235234} + q^2 F_{243523} - q^2 F_{223543} + q^2 F_{324352} + q^2 F_{532432} + (- q^3 - q) F_{234352}\\
    &    + (- q^3 - q) F_{253432} + (- q^3 - q) F_{325432} + (q^4 + q^2 + 1) F_{235432} \\
    A_{14} &= q^2 E_{235234} + q^2 E_{243523} - q^2 E_{223543} + q^2 E_{324352} + q^2 E_{532432} + (- q^3 - q) E_{234352} \\
    &  + (- q^3 - q) E_{253432}   + (- q^3 - q) E_{325432} + (q^4 + q^2 + 1) E_{235432}, \\
    A_{15} &= q^4 F_{23543} - q^3 F_{23435} - q^3 F_{25343} + q^2 F_{32435} + q^2 F_{32534} + q^2 F_{43253} + q^2 F_{53243}\\
    &    - q F_{34325} - q F_{53432} + F_{35432} + (- q^3 - q) F_{32543} \\
    A_{16} &= q^4 E_{34532} - (q^3 - q) E_{35342} - (q^3 - q) E_{43532} + q^2 E_{34235} + q^2 E_{35234} + q^2 E_{43523} + q^2 E_{53423}\\
    &     - q E_{32345} - q E_{43235} - q E_{45323} - q E_{53234} + E_{23453} + (- q^3 - q) E_{34523} \\
    A_{17} &= q^3 F_{2354} - q^2 F_{2435} - q^2 F_{2534} - q^2 F_{3254} + q F_{2543} + q F_{4325} + q F_{5324} - F_{5432} \\
    A_{18} &= q^3 E_{4532} - q^2 E_{4352} - q^2 E_{4523} - q^2 E_{5342} + q E_{3452} + q E_{4235} + q E_{5234} - E_{2345} \\
    A_{19} &= q^6 F_{2354321} + q^4 F_{2343512} + q^4 F_{2534312} + (q^4 - q^2 ) F_{2352341} + (q^4 - q^2 ) F_{2435231}\\
    &    - (q^4 - q^2 ) F_{2235431} + (q^4 - q^2 ) F_{3243521} + (q^4 - q^2 ) F_{5324321} - q^3 F_{2341235} - q^3 F_{2351234} \\
    &    - q^3 F_{2435123} - q^3 F_{2534123} - q^3 F_{3243512} - q^3 F_{3253412} - q^3 F_{4325312} - q^3 F_{5324312} + q^2 F_{2312354}  \\
    &    + q^2 F_{2431235} + q^2 F_{2531234} + q^2 F_{2543123} + q^2 F_{3241235} + q^2 F_{3251234} + q^2 F_{3432512} + q^2 F_{4325123}\\
    &    + q^2 F_{5324123} + q^2 F_{5343212} - q F_{2123543} - q F_{3212354} - q F_{3543212} - q F_{4321235}\\
    &    - q F_{5321234} - q F_{5432123} + (- q^3 - q) F_{3254123} + (q^4 + q^2) F_{2354123} + (q^4 + q^2) F_{3254312} \\
    &    + (- q^5 - q) F_{2354312} + (- q^5 + q) F_{2343521} + (- q^5 + q) F_{2534321} + (- q^5 + q) F_{3254321} + F_{1235432}\\
    A_{20} &= q^6 E_{1235432} - q^5 E_{1234352} - q^5 E_{1253432} - q^5 E_{1325432} + q^4 E_{1235234} + q^4 E_{1243523}\\
    &     + q^4 E_{1324352} + q^4 E_{1532432} - (q^4 - q^2 ) E_{1223543} - q^3 E_{1235423} - q^3 E_{2134235} - q^3 E_{2135234}\\
    &     - q^3 E_{2143523} - q^3 E_{2153423} - q^3 E_{2354213} - q^3 E_{3214352} - q^3 E_{3215342} - q^3 E_{4321532} \\
    &     - q^3 E_{5321432} + q^2 E_{2352134} + q^2 E_{2435213} + q^2 E_{3243521} + q^2 E_{5324321} - q E_{2343521} - q E_{2534321} \\
    &     - q E_{3254321} + (q^3 - q) E_{2123543} + (q^4 + q^2) E_{2134352} + (q^4 + q^2) E_{2135423}\\
    &     + (q^4 + q^2) E_{2153432} + (q^4 + q^2) E_{3215432} + (- q^5 - q) E_{2135432} + E_{2354321} \\
    A_{21} &= - q^3 F_{12341235} + q^3 F_{12354312} - q^3 F_{12534123} + q^3 F_{13212354} - q^3 F_{13253412} + q^3 F_{13543212}  \\ 
    &   + q^3 F_{14321235} - q^3 F_{14325312} + q^3 F_{15321234} + q^3 F_{15432123} + q^3 F_{21235431} - q^3 F_{21352341} \\
    &    - q^3 F_{21435231} - q^3 F_{32143521} - q^3 F_{53214321} - (q^4 + q^2) F_{12235431} - (q^4 + q^2) F_{11235432} \\
    &    + (q^4 + q^2) F_{12352341} + (q^4 + q^2) F_{12354123} + (q^4 + q^2) F_{12435231} + (q^4 + q^2) F_{13243521}  \\
    &     + (q^4 + q^2) F_{15324321} + (q^4 + q^2) F_{21343521} + (q^4 + q^2) F_{21534321} + (q^4 + q^2) F_{32154321} \\
    &    - q (q^2 + 1)^{2} F_{12343521} - q (q^2 + 1)^{2} F_{12534321} - q (q^2 + 1)^{2} F_{13254321} \\
    &    - q (q^4 + q^2 + 1) F_{21354321} + (q^6 + q^4 + q^2 + 1) F_{12354321}
    \end{align*}
and
    \begin{align*}
    A_{22} &= - q^3 E_{12341235} + q^3 E_{12354312}  - q^3 E_{12534123} + q^3 E_{13212354} - q^3 E_{13253412} + q^3 E_{13543212}\\
    &     + q^3 E_{14321235} - q^3 E_{14325312} + q^3 E_{15321234} + q^3 E_{15432123} + q^3 E_{21235431} - q^3 E_{21352341}\\
    &     - q^3 E_{21435231} - q^3 E_{32143521} - q^3 E_{53214321} - (q^4 + q^2) E_{12235431} - (q^4 + q^2) E_{11235432} \\
    &    + (q^4 + q^2) E_{12352341} + (q^4 + q^2) E_{12354123} + (q^4 + q^2) E_{12435231} + (q^4 + q^2) E_{13243521} \\
    &    + (q^4 + q^2) E_{15324321} + (q^4 + q^2) E_{21343521} + (q^4 + q^2) E_{21534321} + (q^4 + q^2) E_{32154321}\\
    &    - q (q^2 + 1)^{2} E_{12343521} - q (q^2 + 1)^{2} E_{12534321} - q (q^2 + 1)^{2} E_{13254321}\\
    &    - q (q^4 + q^2 + 1) E_{21354321} + (q^6 + q^4 + q^2 + 1) E_{12354321}, \\
    A_{23} &= q^4 F_{12354} - q^3 F_{12435} - q^3 F_{12534} - q^3 F_{13254} - q^3 F_{21354} + q^2 F_{12543} + q^2 F_{14325}  + q^2 F_{15324} \\
    &    + q^2 F_{21435} + q^2 F_{21534} + q^2 F_{32154} - q F_{15432} - q F_{21543} - q F_{43215} - q F_{53214} + F_{54321}\\
    A_{24} &= q^4 E_{45321} - q^3 E_{43521} - q^3 E_{45231} - q^3 E_{45312} - q^3 E_{53421} + q^2 E_{34521} + q^2 E_{42351} + q^2 E_{43512}\\
    &     + q^2 E_{45123} + q^2 E_{52341} + q^2 E_{53412} - q E_{23451} - q E_{34512} - q E_{41235} - q E_{51234} + E_{12345} \\
    A_{25} &= q^3 F_{1235} - q^2 F_{1253} - q^2 F_{1325} - q^2 F_{2135} + q F_{1532} + q F_{2153} + q F_{3215} - F_{5321} \\
    A_{26} &= q^3 E_{5321} - q^2 E_{3521} - q^2 E_{5231} - q^2 E_{5312} + q E_{2351} + q E_{3512} + q E_{5123} - E_{1235} \\
    A_{27} &= q^5 F_{123543} - q^4 F_{123435} - q^4 F_{125343} - q^4 F_{213543} + q^3 F_{132435} + q^3 F_{132534} + q^3 F_{143253} + q^3 F_{153243} \\
    &  + q^3 F_{213435} + q^3 F_{215343} - q^2 F_{134325} - q^2 F_{153432} - q^2 F_{321435} - q^2 F_{321534} - q^2 F_{432153}\\
    &  - q^2 F_{532143} + q F_{135432} + q F_{343215} + q F_{534321} + (q^3 + q) F_{321543} + (- q^4 - q^2) F_{132543} - F_{354321}\\
    A_{28} &= q^5 E_{354321} - q^4 E_{354312} - (q^4 - q^2 ) E_{343521} - (q^4 - q^2 ) E_{534321} + q^3 E_{342351} + q^3 E_{352341} + q^3 E_{435231}\\
    & + q^3 E_{534231} - q^2 E_{323541} - q^2 E_{341235} - q^2 E_{351234} - q^2 E_{432351} - q^2 E_{435123} - q^2 E_{532341}\\
    & - q^2 E_{534123} - q^2 E_{543231} + q E_{235431} + q E_{312354} + q E_{431235} + q E_{531234} + q E_{543123} \\
    & + (q^3 - q) E_{343512} + (q^3 - q) E_{534312} + (q^3 + q) E_{354123} + (- q^4 - q^2) E_{354231} - E_{123543} \\
    A_{29} &= q^6 F_{1235432} - q^5 F_{1234352} - q^5 F_{1253432} - q^5 F_{1325432} + q^4 F_{1235234} + q^4 F_{1243523} + q^4 F_{1324352}\\
    &+ q^4 F_{1532432} - (q^4 - q^2 ) F_{1223543} - q^3 F_{1235423} - q^3 F_{2134235} - q^3 F_{2135234} - q^3 F_{2143523} - q^3 F_{2153423} \\
    & - q^3 F_{2354213} - q^3 F_{3214352} - q^3 F_{3215342} - q^3 F_{4321532} - q^3 F_{5321432} + q^2 F_{2352134} + q^2 F_{2435213} \\
    &  + q^2 F_{3243521} + q^2 F_{5324321} - q F_{2343521} - q F_{2534321} - q F_{3254321} + (q^3 - q) F_{2123543} + (q^4 + q^2) F_{2134352} \\
    & + (q^4 + q^2) F_{2135423} + (q^4 + q^2) F_{2153432} + (q^4 + q^2) F_{3215432} + (- q^5 - q) F_{2135432} + F_{2354321} \\
    A_{30} &= q^6 E_{2354321} + q^4 E_{2343512} + q^4 E_{2534312} + (q^4 - q^2 ) E_{2352341} + (q^4 - q^2 ) E_{2435231} - (q^4 - q^2 ) E_{2235431}\\
    & + (q^4 - q^2 ) E_{3243521} + (q^4 - q^2 ) E_{5324321} - q^3 E_{2341235} - q^3 E_{2351234} - q^3 E_{2435123} - q^3 E_{2534123}\\
    & - q^3 E_{3243512} - q^3 E_{3253412} - q^3 E_{4325312} - q^3 E_{5324312} + q^2 E_{2312354} + q^2 E_{2431235} + q^2 E_{2531234}\\
    & + q^2 E_{2543123} + q^2 E_{3241235} + q^2 E_{3251234} + q^2 E_{3432512} + q^2 E_{4325123} + q^2 E_{5324123} + q^2 E_{5343212}\\
    & - q E_{2123543} - q E_{3212354} - q E_{3543212} - q E_{4321235} - q E_{5321234} - q E_{5432123} + (- q^3 - q) E_{3254123}\\
    & + (q^4 + q^2) E_{2354123} + (q^4 + q^2) E_{3254312} + (- q^5 - q) E_{2354312} + (- q^5 + q) E_{2343521} + (- q^5 + q) E_{2534321}\\
    & + (- q^5 + q) E_{3254321} + E_{1235432} \\
   A_{31} &= {q^3} F_{4321} - {q^2} F_{3421} - {q^2} F_{4231} - {q^2} F_{4312} + {q} F_{2341} + {q} F_{3412} + {q} F_{4123} - F_{1234}\\
   A_{32} &= {q^3} E_{1234} - {q^2} E_{1243} - {q^2} E_{1324} - {q^2} E_{2134} + {q} E_{1432} + {q} E_{2143} + {q} E_{3214} - E_{4321}.
\end{align*}    
    This element acts as $q^{10} + q^6 + q^4 + q^2 + 2 + \frac{1}{q^2} + \frac{1}{q^4} + \frac{1}{q^6} + \frac{1}{q^{10}}$ times the identity in the fundamental representation of $\Uq(\so_{10})$.
\end{Theorem}

Using this central element, the method in \cite{CGRS} can be applied to obtain a Markov process. The following theorem states that the resulting process is indeed the Type D ASEP with $n = 5$.

\begin{Theorem}\label{thm:asep}
    Let $C$ denote the central element of $\Uq(\so_{10})$ from Theorem \ref{thm:central-element}. Let $H$ denote the action of $\Delta(C)$ on $\mathbb{C}^{10} \otimes \mathbb{C}^{10}$, so that $H$ is a $100\times 100$ matrix. Define the quantum Hamiltonian $\hat{H} = H - \Lambda \cdot \text{Id}$, where
    $$\Lambda = q^{12} + q^6 + q^4 + q^2 + 2 + \frac{1}{q^2} + \frac{1}{q^4} + \frac{1}{q^6} + \frac{1}{q^{12}}.$$
    There exist linearly independent eigenvectors $g_0, g_1, g_2, g_3$ of $\hat{H}$ with eigenvalue $0$ such that if $G_\delta$ is the diagonal matrix with entries given by $g_\delta$, then removing certain states from $G_\delta^{-1}\hat{H}G_\delta$ results in the generator of the two-site Type D ASEP with parameters $(q, 5, \delta)$, for $\delta = 0, 1, 2, 3$.
\end{Theorem}

In addition, we prove the following generalization of Theorem 3.1 from \cite{Blyschak_2023}:

\begin{Theorem}\label{thm:duality}
    The Type D ASEP with $\delta = 0$ is self-dual with respect to the self-duality function
    $$
    D_{\alpha_1,\alpha_2}^L(\eta,\xi) = D_{\alpha_1}^L(\eta_1,\xi_1)\cdot D_{\alpha_2}^L(\eta_2,\xi_2)
    $$
    where
    $$
    D_{\alpha_i}^L(\xi_i, \eta_i) = \prod_{x=1}^L K_{\eta_{i}^x}\Big(q^{-2\xi_i^x},p_i^x(\xi_i, \eta_i),1,q^2\Big)
    $$
    and
    $$
    p_i^x(\xi_i,\eta_i) = \alpha_i^{-1}q^{-2\left(N_{x-1}^{-}(\xi_i)-N_{x+1}^{+}(\eta_i)\right)+2x-2}.
    $$
\end{Theorem}

Note that here, $\alpha_1$ and $\alpha_2$ are not roots of $\so_{10}$; they are parameters that depend on the reversible measures explained in \cite{Blyschak_2023}.

\section{Proofs}\label{sec:proofs}

The proofs of the first two theorems were assisted by a computer. The code used can be found at \url{https://github.com/e-rohr/Type-D-ASEP}.

\subsection{Proof of Theorem \ref{thm:central-element}}

We find the central element using \cite[Lemma 3.1]{Kuan_2016}, as was done in \cite{kuan2020} for $\Uq(\so_6)$ and $\Uq(\so_8)$. The lemma is restated here for convenience:

\begin{Lemma}
    For each weight $\mu$ of the fundamental representation of a Lie algebra $\mathfrak{g}$, let $v_\mu$ be a vector in the weight space. Suppose $q$ is not a root of unity, and $2\mu$ is always in the root lattice of $\mathfrak{g}$. Let $e_{\mu\lambda}$ and $f_{\lambda\mu}$ be products of $E_i$'s and $F_i$'s in $\Uq(\mathfrak{g})$, respectively, such that $e_{\mu\lambda}$ sends $v_\lambda$ to $v_\mu$ and $f_{\lambda\mu}$ sends $v_\mu$ to $v_\lambda$. If $e_{\mu\lambda}^*$ and $f_{\lambda\mu}^*$ are the corresponding dual elements under $\langle \cdot, \cdot \rangle$, and $\rho$ is half the sum of the positive roots of $\mathfrak{g}$, then 
    \begin{equation}\label{eq:central-element}
        \sum_{\mu\geq\lambda} q^{(\mu - \lambda, \mu)} q^{-(2\rho, \mu)} e_{\mu\lambda}^*q^{H_{-\mu-\lambda}}f_{\lambda\mu}^*
    \end{equation}
    is a central element of $\Uq(\mathfrak{g})$.
    \end{Lemma}

    Several terms in \eqref{eq:central-element} are straightforward to compute; the computations are explained in \cite[Section 2.1.1]{kuan2020}. In particular:
    $$q^{(\mu - \lambda, \mu)} = \begin{cases} q^2 & \lambda = -\mu \\ 1 & \lambda = \mu \\ q & \text{otherwise,} \end{cases} \qquad q^{(-2\rho, \mu)} = \begin{cases} q^{2i - 2n} & \mu = L_i \\ q^{2n - 2i} & \mu = -L_i, \end{cases}$$

$$H_{-2L_i} = H_{n-1} - H_n -2 \sum_{j=i}^{n-1} H_j,$$

and if $i < j$,

$$H_{L_i - L_j} = \sum_{k=i}^{j-1} H_k.$$

We can add or subtract these last two equations to find any $H_{\pm L_i \pm L_j}$.
    
Next, we find $e_{\mu\lambda}$ and $f_{\lambda\mu}$. The following diagram shows the actions of $E_i$ and $F_i$ on $\{v_1, \dots, v_{2n}\}$ for $n = 5$. Recall that the vectors $v_1,\ldots,v_{2n}$ belong to the weight spaces of $L_1,\ldots, L_n,-L_n,\ldots, -L_1$, respectively. 
    
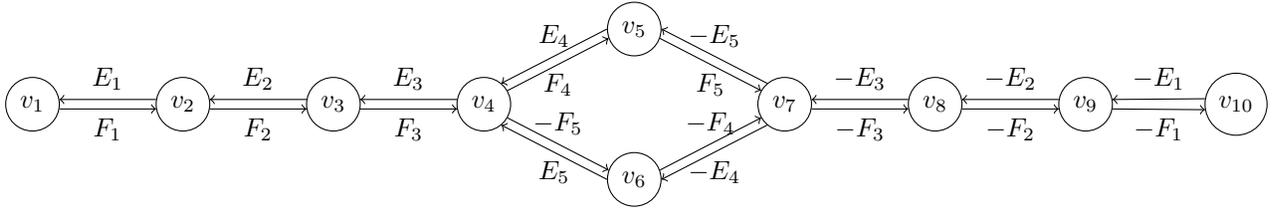
\begin{figure}[ht]
\centering
    \begin{tikzpicture}
    \node[circle, draw] (1) at (-8,0) {$v_{1}$};
    \node[circle, draw] (2) at (-6,0) {$v_{2}$}; 
    \node[circle, draw] (3) at (-4,0) {$v_{3}$}; 
    \node[circle, draw] (4) at (-2,0) {$v_{4}$}; 
    \node[circle, draw] (5) at (0,1) {$v_{5}$}; 
    \node[circle, draw] (-5) at (0,-1) {$v_{6}$}; 
    \node[circle, draw] (-4) at (2,0) {$v_{7}$}; 
    \node[circle, draw] (-3) at (4,0) {$v_{8}$}; 
    \node[circle, draw] (-2) at (6,0) {$v_{9}$}; 
    \node[circle, draw] (-1) at (8,0) {$v_{10}$}; 

    \draw (2.170) edge [->]  node [above] {$E_1$} (1.10);
    \draw (3.170) edge [->]  node [above] {$E_2$} (2.10);
    \draw (4.170) edge [->]  node [above] {$E_3$} (3.10);
    \draw (5.180) edge [->]  node [above] {$E_4$} (4.50);
    \draw (-5.180) edge [->]  node [below] {$E_5$} (4.310);
    \draw (-4.130) edge [->]  node [above] {$-E_5$} (5.0);
    \draw (-4.230) edge [->]  node [below] {$-E_4$} (-5.0);
    \draw (-3.170) edge [->]  node [above] {$-E_3$} (-4.10);
    \draw (-2.170) edge [->]  node [above] {$-E_2$} (-3.10);
    \draw (-1.170) edge [->]  node [above] {$-E_1$} (-2.10);

    \draw (1.350) edge [->]  node [below] {$F_1$} (2.190);
    \draw (2.350) edge [->]  node [below] {$F_2$} (3.190);
    \draw (3.350) edge [->]  node [below] {$F_3$} (4.190);
    \draw (4.30) edge [->]  node [below] {$F_4$} (5.200);
    \draw (4.330) edge [->]  node [above] {$-F_5$} (-5.160);
    \draw (5.340) edge [->]  node [below] {$F_5$} (-4.150);
    \draw (-5.20) edge [->]  node [above] {$-F_4$} (-4.210);
    \draw (-4.350) edge [->]  node [below] {$-F_3$} (-3.190);
    \draw (-3.350) edge [->]  node [below] {$-F_2$} (-2.190);
    \draw (-2.350) edge [->]  node [below] {$-F_1$} (-1.190);
    \end{tikzpicture}
\end{figure}

For example, for $\mu = L_3$ and $\lambda = -L_2$, we have $e_{\mu\lambda} = - E_3 E_4 E_5 E_3 E_2 $ and $f_{\lambda \mu} = F_2 F_3 F_5 F_4 F_3$. Note that there is no difference between taking the top and bottom path because $E_4, E_5$ commute and $F_4, F_5$ commute.

The procedure for obtaining the dual elements $e_{\mu\lambda}^*$ and $f_{\lambda \mu}^*$ is described in depth in \cite{kuan2020} and summarized here. Given $e_{\mu\lambda}$ of the form $E_{x_1} \dots E_{x_k}$, create the set of elements $\{E_{\sigma(x_1)} \dots E_{\sigma(x_k)}\}_\sigma$, where $\sigma$ ranges through all possible permutations. Some of these elements may be linearly dependent due to the relations in the quantum group. Thus, take a basis $e_1 = e_{\mu\lambda}, e_2, \dots, e_m$ of the span of these elements, and create a corresponding basis $f_1, \dots, f_m$, where each $f_i$ is the same as $e_i$ but with $F_{j}$ replacing each $E_{j}$. Form a matrix $M$ of $q$-pairings such that $M_{ij} = \langle e_i, f_j \rangle$. Then $e_{\mu\lambda}^*$ is the dot product of the first row of $M^{-1}$ with the column vector $(f_1, \dots, f_m)$. Obtaining $f_{\lambda \mu}^*$ is similar.

The calculation for $\Uq(\so_{10})$ is done with the aid of a computer. The same code could be used for calculating a central element of $\Uq(\so_{2n})$ for larger $n$. For $\Uq(\so_6)$ and $\Uq(\so_8)$, the matrices are small enough to perform computations symbolically, but computing determinants and inverses of large symbolic matrices is intractable. These difficulties are especially prevalent in the computation of the dual elements $e_{\mu\lambda}^*$ and $f_{\lambda \mu}^*$. As such, for $\Uq(\so_{10})$, the process of obtaining the dual elements is done with a numerical value of $q = 10$ to speed up computation. 

\subsection{Proof of Theorem \ref{thm:asep}}

Under the appropriate basis, the quantum Hamiltonian $\hat{H}$ decomposes as a direct sum of one $10 \times 10$ block, forty $2 \times 2$ blocks, and ten $1 \times 1$ blocks. 

The $1 \times 1$ blocks all have entry zero, so taking four of these blocks corresponds to the four $1 \times 1$ blocks from the generator matrix $\mathcal{L}$ in \eqref{eq:generator}.

The $2 \times 2$ blocks are
$$\begin{bmatrix}
    - q^{10} + 2 q^{8} - q^{6} - \frac{1}{q^{8}} + \frac{2}{q^{10}} - \frac{1}{q^{12}} & \frac{\left(q^{2} - 1\right)^{2} \left(q^{18} + 1\right)}{q^{11}}\\\frac{\left(q^{2} - 1\right)^{2} \left(q^{18} + 1\right)}{q^{11}} & - q^{12} + 2 q^{10} - q^{8} - \frac{1}{q^{6}} + \frac{2}{q^{8}} - \frac{1}{q^{10}}
\end{bmatrix},$$
which have eigenvector $\begin{pmatrix} q \\ 1\end{pmatrix}$, and conjugating by the corresponding diagonal matrix $\begin{pmatrix} q & 0\\ 0 & 1\end{pmatrix}$ results in one of the $2 \times 2$ blocks $\mathcal{L}_2$ in the generator \eqref{eq:generator} multiplied by a constant factor of $r^2$.

Let
\begin{align*}
    B_1 &= - q^{5} + 3 q^{3} - 3 q + \tfrac{1}{q} - \tfrac{2}{q^{5}} + \tfrac{4}{q^{7}} - \tfrac{2}{q^{9}}, \\
    B_2 &= - 2 q^{3} + 4 q - \tfrac{2}{q} + \tfrac{1}{q^{5}} - \tfrac{3}{q^{7}} + \tfrac{3}{q^{9}} - \tfrac{1}{q^{11}}, \\
    B_3 &= q^{10} - 2 q^{8} + q^{6} - 2 q^{2} + 4 - \tfrac{2}{q^{2}} + \tfrac{1}{q^{6}} - \tfrac{2}{q^{8}} + \tfrac{1}{q^{10}}.
\end{align*}
Then the $10 \times 10$ block has form $U^T + D + U$, where
$$U = \left[\begin{matrix}0 & B_1 & qB_1 & q^2 B_1 & q^3 B_1 & B_3 & q^6 B_1 & q^5 B_1 & q^4 B_1 & q^3 B_1 \\
0 & 0 & q^2 B_1 & q^3 B_1 & q^4 B_1 & B_2 & B_3 & q^6 B_1 & q^5 B_1 & q^4 B_1\\
0 & 0 & 0 & q^4 B_1 & q^5 B_1 & q B_2 & B_2 & B_3 & q^6 B_1 & q^5 B_1 \\
0 & 0 & 0 & 0 & q^6 B_1 & q^2 B_3 & q B_2 & B_2 & B_3 & q^6 B_1 \\
0 & 0 & 0 & 0 & 0 & q^3 B_3 & q^2 B_2 & q B_2 & B_2 & B_3 \\
0 & 0 & 0 & 0 & 0 & 0 & q^6 B_2 & q^5 B_2 & q^4 B_2 & q^3 B_2 \\
0 & 0 & 0 & 0 & 0 & 0 & 0 & q^4 B_2 & q^3 B_2 & q^2 B_2 \\
0 & 0 & 0 & 0 & 0 & 0 & 0 & 0 & q^2 B_2 & q B_2 \\0 & 0 & 0 & 0 & 0 & 0 & 0 & 0 & 0 & B_2 \\0 & 0 & 0 & 0 & 0 & 0 & 0 & 0 & 0 & 0\end{matrix}\right],$$
and $D$ is the diagonal matrix with entries
\begin{footnotesize}
\begin{equation*}
\begin{split}
\Bigg\{- q^{10} + 2 q^{8} - q^{6} - q^{4} + 3 q^{2} - 3 + \frac{1}{q^{2}} - \frac{2}{q^{6}} + \frac{3}{q^{8}} - \frac{1}{q^{12}},
- q^{10} + 2 q^{8} - 2 q^{6} + 3 q^{4} - 3 q^{2} + 1 - \frac{2}{q^{4}} + \frac{4}{q^{6}} - \frac{3}{q^{8}} + \frac{2}{q^{10}} - \frac{1}{q^{12}},\\
- q^{10} + q^{8} + 2 q^{6} - 3 q^{4} + q^{2} - \frac{2}{q^{2}} + \frac{4}{q^{4}} - \frac{2}{q^{6}} - \frac{1}{q^{8}} + \frac{2}{q^{10}} - \frac{1}{q^{12}},
- 2 q^{10} + 5 q^{8} - 4 q^{6} + q^{4} - 2 + \frac{4}{q^{2}} - \frac{2}{q^{4}} - \frac{1}{q^{8}} + \frac{2}{q^{10}} - \frac{1}{q^{12}},\\
- q^{12} + 2 q^{10} - q^{8} - 2 q^{2} + 4 - \frac{2}{q^{2}} - \frac{1}{q^{8}} + \frac{2}{q^{10}} - \frac{1}{q^{12}},
- q^{12} + 3 q^{8} - 2 q^{6} + q^{2} - 3 + \frac{3}{q^{2}} - \frac{1}{q^{4}} - \frac{1}{q^{6}} + \frac{2}{q^{8}} - \frac{1}{q^{10}},\\
- q^{12} + 2 q^{10} - 3 q^{8} + 4 q^{6} - 2 q^{4} + 1 - \frac{3}{q^{2}} + \frac{3}{q^{4}} - \frac{2}{q^{6}} + \frac{2}{q^{8}} - \frac{1}{q^{10}},
- q^{12} + 2 q^{10} - q^{8} - 2 q^{6} + 4 q^{4} - 2 q^{2} + \frac{1}{q^{2}} - \frac{3}{q^{4}} + \frac{2}{q^{6}} + \frac{1}{q^{8}} - \frac{1}{q^{10}},\\
- q^{12} + 2 q^{10} - q^{8} - 2 q^{4} + 4 q^{2} - 2 + \frac{1}{q^{4}} - \frac{4}{q^{6}} + \frac{5}{q^{8}} - \frac{2}{q^{10}},
- q^{12} + 2 q^{10} - q^{8} - 2 q^{2} + 4 - \frac{2}{q^{2}} - \frac{1}{q^{8}} + \frac{2}{q^{10}} - \frac{1}{q^{12}}
\Bigg\}.
\end{split}
\end{equation*}
\end{footnotesize}

The $10\times 10$ block has rank $6$ and has four linearly independent eigenvectors 
$$
g_0=
\left(
\begin{array}{c}
0\\
0\\
0\\
-q^2\\
q\\
0\\
0\\
0\\
-1\\
q
\end{array}
\right), \quad
g_1=\left(
\begin{array}{c}
0\\
0\\
-q^2\\
q\\
0\\
0\\
0\\
-1\\
q\\
0\\
\end{array}
\right), \quad
g_2=
\left(
\begin{array}{c}
0\\
-q^2\\
q\\
0\\
0\\
0\\
-1\\
q\\
0\\
0
\end{array}
\right), \quad
g_3=\left(
\begin{array}{c}
-q^2\\
q\\
0\\
0\\
0\\
-1\\
q\\
0\\
0\\
0
\end{array}
\right).
$$
These are the four choices of ground state vector corresponding to the Type D ASEPs with parameters $(q, 5, 0), (q, 5, 1), (q, 5, 2)$, and $(q, 5, 3)$, respectively. 

For each $\delta = 0, 1, 2, 3$, set $L_\delta = G_\delta^{-1} \hat{H} G_\delta$. Then,
$$
L_0 = \left[\begin{matrix}* & B_{1} & q B_{1} & \infty & \infty & B_{3} & q^{6} B_{1} & q^{5} B_{1} & \infty & \infty\\B_{1} & * & q^{2} B_{1} & \infty & \infty & B_{2} & B_{3} & q^{6} B_{1} & \infty & \infty\\q B_{1} & q^{2} B_{1} & * & \infty & \infty & q B_{2} & B_{2} & B_{3} & \infty & \infty\\0 & 0 & 0 & * & - q^{5} B_{1} & 0 & 0 & 0 & \frac{B_{3}}{q^{2}} & - q^{5} B_{1}\\0 & 0 & 0 & - q^{7} B_{1} & * & 0 & 0 & 0 & - \frac{B_{2}}{q} & B_{3}\\B_{3} & B_{2} & q B_{2} & \infty & \infty & * & q^{6} B_{2} & q^{5} B_{2} & \infty & \infty\\q^{6} B_{1} & B_{3} & B_{2} & \infty & \infty & q^{6} B_{2} & * & q^{4} B_{2} & \infty & \infty\\q^{5} B_{1} & q^{6} B_{1} & B_{3} & \infty & \infty & q^{5} B_{2} & q^{4} B_{2} & * & \infty & \infty\\0 & 0 & 0 & q^{2} B_{3} & - q B_{2} & 0 & 0 & 0 & * & - q B_{2}\\0 & 0 & 0 & - q^{7} B_{1} & B_{3} & 0 & 0 & 0 & - \frac{B_{2}}{q} & *\end{matrix}\right],
$$

$$
L_1 = \left[\begin{matrix}* & B_{1} & \infty & \infty & q^{3} B_{1} & B_{3} & q^{6} B_{1} & \infty & \infty & q^{3} B_{1}\\B_{1} & * & \infty & \infty & q^{4} B_{1} & B_{2} & B_{3} & \infty & \infty & q^{4} B_{1}\\0 & 0 & * & - q^{3} B_{1} & 0 & 0 & 0 & \frac{B_{3}}{q^{2}} & - q^{5} B_{1} & 0\\0 & 0 & - q^{5} B_{1} & * & 0 & 0 & 0 & - \frac{B_{2}}{q} & B_{3} & 0\\q^{3} B_{1} & q^{4} B_{1} & \infty & \infty & * & q^{3} B_{2} & q^{2} B_{2} & \infty & \infty & B_{3}\\B_{3} & B_{2} & \infty & \infty & q^{3} B_{2} & * & q^{6} B_{2} & \infty & \infty & q^{3} B_{2}\\q^{6} B_{1} & B_{3} & \infty & \infty & q^{2} B_{2} & q^{6} B_{2} & * & \infty & \infty & q^{2} B_{2}\\0 & 0 & q^{2} B_{3} & - q B_{2} & 0 & 0 & 0 & * & - q^{3} B_{2} & 0\\0 & 0 & - q^{7} B_{1} & B_{3} & 0 & 0 & 0 & - q B_{2} & * & 0\\q^{3} B_{1} & q^{4} B_{1} & \infty & \infty & B_{3} & q^{3} B_{2} & q^{2} B_{2} & \infty & \infty & *\end{matrix}\right],
$$

$$
L_2 = \left[\begin{matrix}* & \infty & \infty & q^{2} B_{1} & q^{3} B_{1} & B_{3} & \infty & \infty & q^{4} B_{1} & q^{3} B_{1}\\0 & * & - q B_{1} & 0 & 0 & 0 & \frac{B_{3}}{q^{2}} & - q^{5} B_{1} & 0 & 0\\0 & - q^{3} B_{1} & * & 0 & 0 & 0 & - \frac{B_{2}}{q} & B_{3} & 0 & 0\\q^{2} B_{1} & \infty & \infty & * & q^{6} B_{1} & q^{2} B_{2} & \infty & \infty & B_{3} & q^{6} B_{1}\\q^{3} B_{1} & \infty & \infty & q^{6} B_{1} & * & q^{3} B_{2} & \infty & \infty & B_{2} & B_{3}\\B_{3} & \infty & \infty & q^{2} B_{2} & q^{3} B_{2} & * & \infty & \infty & q^{4} B_{2} & q^{3} B_{2}\\0 & q^{2} B_{3} & - q B_{2} & 0 & 0 & 0 & * & - q^{5} B_{2} & 0 & 0\\0 & - q^{7} B_{1} & B_{3} & 0 & 0 & 0 & - q^{3} B_{2} & * & 0 & 0\\q^{4} B_{1} & \infty & \infty & B_{3} & B_{2} & q^{4} B_{2} & \infty & \infty & * & B_{2}\\q^{3} B_{1} & \infty & \infty & q^{6} B_{1} & B_{3} & q^{3} B_{2} & \infty & \infty & B_{2} & *\end{matrix}\right],
$$
and
$$
L_3 = \left[\begin{matrix}* & - \frac{B_{1}}{q} & 0 & 0 & 0 & \frac{B_{3}}{q^{2}} & - q^{5} B_{1} & 0 & 0 & 0\\- q B_{1} & * & 0 & 0 & 0 & - \frac{B_{2}}{q} & B_{3} & 0 & 0 & 0\\\infty & \infty & * & q^{4} B_{1} & q^{5} B_{1} & \infty & \infty & B_{3} & q^{6} B_{1} & q^{5} B_{1}\\\infty & \infty & q^{4} B_{1} & * & q^{6} B_{1} & \infty & \infty & B_{2} & B_{3} & q^{6} B_{1}\\\infty & \infty & q^{5} B_{1} & q^{6} B_{1} & * & \infty & \infty & q B_{2} & B_{2} & B_{3}\\q^{2} B_{3} & - q B_{2} & 0 & 0 & 0 & * & - q^{7} B_{2} & 0 & 0 & 0\\- q^{7} B_{1} & B_{3} & 0 & 0 & 0 & - q^{5} B_{2} & * & 0 & 0 & 0\\\infty & \infty & B_{3} & B_{2} & q B_{2} & \infty & \infty & * & q^{2} B_{2} & q B_{2}\\\infty & \infty & q^{6} B_{1} & B_{3} & B_{2} & \infty & \infty & q^{2} B_{2} & * & B_{2}\\\infty & \infty & q^{5} B_{1} & q^{6} B_{1} & B_{3} & \infty & \infty & q B_{2} & B_{2} & *\end{matrix}\right],
$$
where the diagonal entries * are just the entries of $D$. Remove any rows containing $\infty$ and any columns containing $0$ from $L_\delta$, and divide all entries by $r^2$ to obtain a $4 \times 4$ matrix $\tilde{L}_\delta$. These are

\setlength{\arraycolsep}{2pt}
\renewcommand*{\arraystretch}{2.5}
$$
\tilde{L}_0 = \left[\begin{matrix}\frac{- 2 q^{18} + q^{16} - 2 q^{8} - 1}{q^{10}} & \frac{q^{10} - q^{8} + 2}{q^{2}} & \frac{q^{16} - 2 q^{8} + 1}{q^{10}} & \frac{q^{10} - q^{8} + 2}{q^{2}}\\
q^{10} - q^{8} + 2 & \frac{- q^{20} - 2q^{10} - 1}{q^{10}} & \frac{2q^{10} - q^{2} + 1}{q^{10}} & \frac{q^{16} - 2q^{8} + 1}{q^{8}}\\
\frac{q^{16} - 2 q^{8} + 1}{q^{6}} & \frac{2 q^{10} - q^{2} + 1}{q^{8}} & \frac{- q^{18} - 2 q^{10} + q^{2} - 2}{q^{8}} & \frac{2 q^{10} - q^{2} + 1}{q^{8}}\\
q^{10} - q^{8} + 2 & \frac{q^{16} - 2q^{8} + 1}{q^{8}} & \frac{2q^{10} - q^{2} + 1}{q^{10}} &  \frac{- q^{20} - 2q^{10} - 1}{q^{10}}\end{matrix}\right],
$$

\setlength{\arraycolsep}{5pt}
\renewcommand*{\arraystretch}{2}
$$
\tilde{L}_1 = \left[\begin{matrix}\frac{- q^{18} - q^{16} + q^{14} - 2 q^{6} - 1}{q^{10}} & \frac{q^{10} - q^{8} + 2}{q^{4}} & \frac{q^{16} - 2 q^{8} + 1}{q^{10}} & \frac{q^{10} - q^{8} + 2}{q^{2}}\\
\frac{q^{10} - q^{8} + 2}{q^{2}} & \frac{- 2 q^{18} + q^{16} - 2 q^{8} - 1}{q^{10}} & \frac{2q^{10} - {q^{2}} + 1}{q^{10}} & \frac{q^{16} - 2q^8 + 1}{q^{8}}\\
\frac{q^{16} - 2 q^{8} + 1}{q^{6}} & \frac{2 q^{10} - q^{2} + 1}{q^{8}} & \frac{- q^{18} - 2 q^{12} + q^{4} - q^{2} - 1}{q^{8}} & \frac{2 q^{10} - q^{2} + 1}{q^{6}}\\
q^{10} - q^{8} + 2 & \frac{q^{16} - 2q^8 + 1}{q^{8}} & \frac{2 q^{10} - q^{2} + 1}{q^{8}} & \frac{- q^{18} - 2 q^{10} + q^{2} - 2}{q^{8}}\end{matrix}\right],
$$

\setlength{\arraycolsep}{.1pt}
\renewcommand*{\arraystretch}{2}
$$
\tilde{L}_2 = \left[\begin{matrix}\frac{- q^{18} - q^{14} + q^{12} - 2 q^{4} - 1}{q^{10}} & \frac{q^{10} - q^{8} + 2}{q^{6}} & \frac{q^{16} - 2 q^{8} + 1}{q^{10}} & \frac{q^{10} - q^{8} + 2}{q^{2}}\\
\frac{q^{10} - q^{8} + 2}{q^{4}} & \frac{- q^{18} - q^{16} + q^{14} - 2 q^{6} - 1}{q^{10}} & \frac{2q^{10} - {q^{2}} + 1}{q^{10}} & \frac{q^{16} - 2q^8 + 1}{q^{8}}\\
\frac{q^{16} - 2 q^{8} + 1}{q^{6}} & \frac{2 q^{10} - q^{2} + 1}{q^{8}} & \frac{- q^{18} - 2 q^{14} + q^{6} - q^{4} - 1}{q^{8}} & \frac{2 q^{10} - q^{2} + 1}{q^{4}}\\q^{10} - q^{8} + 2 & \frac{q^{16} - 2q^8 + 1}{q^{8}} & \frac{2 q^{10} - q^{2} + 1}{q^{6}} & \frac{- q^{18} - 2 q^{12} + q^{4} - q^{2} - 1}{q^{8}}\end{matrix}\right],
$$
and
\setlength{\arraycolsep}{.1pt}
\renewcommand*{\arraystretch}{2}
$$
\tilde{L}_3 = \left[\begin{matrix}\frac{- q^{18} - q^{12} + q^{10} - 2q^{2} - 1}{q^{10}} & \frac{q^{10} - q^{8} + 2}{q^{8}} & \frac{q^{16} - 2 q^{8} + 1}{q^{10}} & \frac{q^{10} - q^{8} + 2}{q^{2}}\\
\frac{q^{10} - q^{8} + 2}{q^{6}} & \frac{- q^{18} - q^{14} + q^{12} - 2 q^{4} - 1}{q^{10}} & \frac{2q^{10} - {q^{2}} + 1}{q^{10}} & \frac{q^{16} - 2q^8 + 1}{q^{8}}\\
\frac{q^{16} - 2 q^{8} + 1}{q^{6}} & \frac{2 q^{10} - q^{2} + 1}{q^{8}} & \frac{- q^{18} - 2 q^{16} + q^{8} - q^{6} - 1}{q^{8}} & \frac{2 q^{10} - q^{2} + 1}{q^{2}}\\
q^{10} - q^{8} + 2 & \frac{q^{16} - 2q^8 + 1}{q^{8}} & \frac{2 q^{10} - q^{2} + 1}{q^{4}} & \frac{- q^{18} - 2 q^{14} + q^{6} - q^{4} - 1}{q^{8}}\end{matrix}\right].
$$
These matrices are written with respect to the ordered bases $(v_4 \otimes v_9, v_5 \otimes v_{10}, v_9 \otimes v_4, v_{10} \otimes v_5), (v_3 \otimes v_8, v_4 \otimes v_9, v_8 \otimes v_3, v_9 \otimes v_4), (v_2 \otimes v_7, v_3 \otimes v_8, v_2 \otimes v_7, v_8 \otimes v_3), (v_1 \otimes v_6, v_2 \otimes v_7, v_6 \otimes v_1, v_2 \otimes v_7)$, respectively.

Observe that $\tilde{L}_\delta$ is exactly the $4 \times 4$ block $\mathcal{L}_1$ in the generator \eqref{eq:generator} with parameter $(q, 5, \delta)$, as desired.

\subsection{Proof of Theorem \ref{thm:duality}}

This theorem generalizes Theorem 3.1 from \cite{Blyschak_2023} for all $n$. The inductive step remains the same as from that paper, so it remains to verify the base case of $L = 2$ for general $n$. 

Let $\mathcal{L}$ be the $16 \times 16$ generator matrix for the two-site model with parameters $(q, n, 0)$ from \eqref{eq:generator}. Let $\mathcal{D}$ be the $16 \times 16$ matrix whose rows and columns are indexed by the states of Type D ASEP (in the same order as for $\mathcal{L}$) and whose $(\eta, \xi)$ entry is $D^L_{\alpha_1, \alpha_2}(\eta, \xi)$. \\
\newpage
Setting $D_i^j = 1-\tfrac{q^j}{\alpha_i}$, we have

\setlength{\arraycolsep}{2.5pt}
\renewcommand*{\arraystretch}{2}

$$
\D = \begin{scriptsize}
    \left[\begin{array}{cccccccccccccccc}D_1^2 D_2^2 & D_2^2 & 1 & D_1^2 & D_1^2 & 1 & D_2^2 & 1 & D_1^2 D_2^2 & D_1^2 & D_1^2 D_2^2 & D_2^2 & 1 & D_1^2 & D_2^2 & D_1^2 D_2^2\\D_2^2 & D_1^4 D_2^2 & D_1^4 & 1 & 1 & D_1^4 & D_2^2 & 1 & D_1^2 D_2^2 & D_1^2 & D_2^2 & D_1^4 D_2^2 & 1 & D_1^2 & D_2^2 & D_1^2 D_2^2\\1 & D_1^4 & D_1^4 D_2^4 & D_2^4 & 1 & D_1^4 & 1 & D_2^4 & D_1^2 & D_1^2 D_2^4 & D_2^2 & D_1^4 D_2^2 & 1 & D_1^2 & D_2^2 & D_1^2 D_2^2\\D_1^2 & 1 & D_2^4 & D_1^2 D_2^4 & D_1^2 & 1 & 1 & D_2^4 & D_1^2 & D_1^2 D_2^4 & D_1^2 D_2^2 & D_2^2 & 1 & D_1^2 & D_2^2 & D_1^2 D_2^2\\D_1^2 & 1 & 1 & D_1^2 & D_1^2 & 1 & 1 & 1 & D_1^2 & D_1^2 & D_1^2 & 1 & 1 & D_1^2 & 1 & D_1^2\\1 & D_1^4 & D_1^4 & 1 & 1 & D_1^4 & 1 & 1 & D_1^2 & D_1^2 & 1 & D_1^4 & 1 & D_1^2 & 1 & D_1^2\\D_2^2 & D_2^2 & 1 & 1 & 1 & 1 & D_2^2 & 1 & D_2^2 & 1 & D_2^2 & D_2^2 & 1 & 1 & D_2^2 & D_2^2\\1 & 1 & D_2^4 & D_2^4 & 1 & 1 & 1 & D_2^4 & 1 & D_2^4 & D_2^2 & D_2^2 & 1 & 1 & D_2^2 & D_2^2\\D_1^4 D_2^2 & D_1^4 D_2^2 & D_1^4 & D_1^4 & D_1^4 & D_1^4 & D_2^2 & 1 & D_1^2 D_1^4 D_2^2 & D_1^2 D_1^4 & D_1^4 D_2^2 & D_1^4 D_2^2 & 1 & D_1^2 D_1^4 & D_2^2 & D_1^2 D_1^4 D_2^2\\D_1^4 & D_1^4 & D_1^4 D_2^4 & D_1^4 D_2^4 & D_1^4 & D_1^4 & 1 & D_2^4 & D_1^2 D_1^4 & D_1^2 D_1^4 D_2^4 & D_1^4 D_2^2 & D_1^4 D_2^2 & 1 & D_1^2 D_1^4 & D_2^2 & D_1^2 D_1^4 D_2^2\\D_1^2 D_2^4 & D_2^4 & D_2^4 & D_1^2 D_2^4 & D_1^2 & 1 & D_2^4 & D_2^4 & D_1^2 D_2^4 & D_1^2 D_2^4 & D_1^2 D_2^2 D_2^4 & D_2^2 D_2^4 & 1 & D_1^2 & D_2^2 D_2^4 & D_1^2 D_2^2 D_2^4\\D_2^4 & D_1^4 D_2^4 & D_1^4 D_2^4 & D_2^4 & 1 & D_1^4 & D_2^4 & D_2^4 & D_1^2 D_2^4 & D_1^2 D_2^4 & D_2^2 D_2^4 & D_1^4 D_2^2 D_2^4 & 1 & D_1^2 & D_2^2 D_2^4 & D_1^2 D_2^2 D_2^4\\1 & 1 & 1 & 1 & 1 & 1 & 1 & 1 & 1 & 1 & 1 & 1 & 1 & 1 & 1 & 1\\D_1^4 & D_1^4 & D_1^4 & D_1^4 & D_1^4 & D_1^4 & 1 & 1 & D_1^2 D_1^4 & D_1^2 D_1^4 & D_1^4 & D_1^4 & 1 & D_1^2 D_1^4 & 1 & D_1^2 D_1^4\\D_2^4 & D_2^4 & D_2^4 & D_2^4 & 1 & 1 & D_2^4 & D_2^4 & D_2^4 & D_2^4 & D_2^2 D_2^4 & D_2^2 D_2^4 & 1 & 1 & D_2^2 D_2^4 & D_2^2 D_2^4\\D_1^4 D_2^4 & D_1^4 D_2^4 & D_1^4 D_2^4 & D_1^4 D_2^4 & D_1^4 & D_1^4 & D_2^4 & D_2^4 & D_1^2 D_1^4 D_2^4 & D_1^2 D_1^4 D_2^4 & D_1^4 D_2^2 D_2^4 & D_1^4 D_2^2 D_2^4 & 1 & D_1^2 D_1^4 & D_2^2 D_2^4 & D_1^2 D_1^4 D_2^2 D_2^4\end{array}\right]
\end{scriptsize}.
$$

One can verify the duality relation
$$
\mathcal{L}\D = \D\mathcal{L}^T
$$
for the base case of $L = 2$, completing the proof.

\bibliographystyle{alpha}
\bibliography{references}

\end{document}